\documentclass{amsart}
\usepackage{amssymb, pst-node, pstricks}

\newtheorem{theorem}{Theorem}[section]
\newtheorem{proposition}[theorem]{Proposition}
\newtheorem{corollary}[theorem]{Corollary}
\newtheorem{definition}[theorem]{Definition}
\newtheorem{conjecture}[theorem]{Conjecture}
\newtheorem{remark}[theorem]{Remark}
\newtheorem{lemma}[theorem]{Lemma}

\psset{unit=1pt,arrowsize=4pt,linewidth=.5pt,linecolor=black}
\def\fns{\footnotesize}

\def\wnot{w_{\rm o}}

\def\mytilde{\kern-.015in\hbox{\lower.03in\hbox{\~{}}}\kern-.01in}
\def\mystrut{\vbox to .15in{}}
\newcommand{\roots}{\Phi}
\newcommand{\subroots}{\Delta}
\newcommand{\posroots}{\Phi_+}
\newcommand{\negroots}{\Phi_-}
\newcommand{\inv}[1]{I_{\roots}(#1)}

\newcommand{\fl}{\mbox{\rm fl}}
\newcommand{\bases}{\mathcal{B}}
\newcommand{\bad}{\mathrm{BAD}}
\thicklines

\begin{document}

\title[Smoothness of Schubert varieties via patterns]{Smoothness of Schubert varieties \\ via patterns in root subsystems}

\author{Sara Billey \and Alexander Postnikov}

\address{Department of Mathematics, University of Washington, Seattle, WA}
\email{billey@math.washington.edu}
\address{Department of Mathematics, M.I.T., Cambridge, MA 02139}
\email{apost@math.mit.edu}


\date{July 28, 2003}

\thanks{S.B.\ was supported by NSF grant number DMS-9983797.}
\thanks{A.P.\ was supported by NSF grant number DMS-0201494.} 

\begin{abstract}
The aim of this article is to present a smoothness criterion for
Schubert varieties in generalized flag manifolds $G/B$ in terms of
patterns in root systems.  We generalize Lakshmibai-Sandhya's
well-known result that says that a Schubert variety in $SL(n)/B$ is
smooth if and only if the corresponding permutation avoids the
patterns $3412$ and $4231$.  Our criterion is formulated uniformly in
general Lie theoretic terms.  We define a notion of pattern in Weyl
group elements and show that a Schubert variety is smooth (or
rationally smooth) if and only if the corresponding element of the
Weyl group avoids a certain finite list of patterns.  These forbidden
patterns live only in root subsystems with star-shaped Dynkin
diagrams.  In the simply-laced case the list of forbidden patterns is
especially simple: besides two patterns of type $A_3$ that appear in
Lakshmibai-Sandhya's criterion we only need one additional forbidden
pattern of type $D_4$.  In terms of these patterns, the only
difference between smoothness and rational smoothness is a single
pattern in type $B_{2}$.  Remarkably, several other important classes
of elements in Weyl groups can also be described in terms of forbidden
patterns.  For example, the fully commutative elements in Weyl groups
have such a characterization.  In order to prove our criterion we used
several known results for the classical types.  For the exceptional
types, our proof is based on computer verifications.  In order to
conduct such a verification for the computationally challenging type
$E_8$, we derived several general results on Poincar\'e polynomials of
cohomology rings of Schubert varieties based on parabolic
decomposition, which have an independent interest.
\end{abstract}

\maketitle

\section{Introduction}

Let $G$ be a semisimple simply-connected complex Lie group
and $B$ be a Borel subgroup.   The generalized flag manifold
$G/B$ decomposes into a disjoint union of {\it Schubert cells\/}
$BwB/B$, labeled by elements $w$ of the corresponding Weyl
group $W$.  The {\it Schubert varieties\/} $X_w=\overline{BwB/B}$ 
are the closures of the Schubert cells.  
A classical question of Schubert calculus is:
{\it For which elements $w$ in the Weyl group $W$, is the Schubert 
variety $X_w$ smooth?}

This question has a particularly nice answer for $G=SL(n)$.  In this 
case the Weyl group is the {symmetric group\/} $W=S_n$
of permutations of $n$ letters.
For a permutation $w=w_1\,w_2\cdots w_n$ in $S_n$ and another
permutation $\sigma=\sigma_1\,\sigma_2\cdots\sigma_k$ in $S_k$, 
with $k\leq n$, 
we say that $w$ {\it contains the pattern\/} $\sigma$ if there is 
a sequence $1\leq p_1<\dots <p_k\leq n$ such that
$w_{p_i}> w_{p_j}$ if and only if $\sigma_i > \sigma_j$ for all 
$1\leq i<j \leq k$.  In other words, $w$ contains the pattern $\sigma$
if there is a subsequence in $w$ of size $k$ with the same relative
order of elements as in $\sigma$.  If $w$ does not contain the pattern
$\sigma$, then we say that $w$ {\it avoids the pattern\/} $\sigma$.

\begin{theorem} \textrm{(Lakshmibai-Sandhya~\cite{LSa})}
For a permutation
$w\in S_n$, the Schubert variety $X_w$ in $SL(n)/B$ is smooth
if and only if $w$ avoids the patterns $3412$ and $4231$.
\end{theorem}
There are several general approaches to determining smoothness of Schubert
varieties.  See Billey and Lakshmibai~\cite{BL} for a survey of 
known results. 
Kazhdan and Lusztig defined a weaker condition called {\it rational
smoothness}.  Rational smoothness can be interpreted in terms of
Kazhdan-Lusztig polynomials~\cite{KL1},~\cite{KL2}. A Schubert variety
is rationally smooth whenever certain Kazhdan-Lusztig polynomials are
trivial.  Kumar~\cite{Kum} presented smoothness and rational
smoothness criteria in terms of the nil Hecke ring, defined
in~\cite{KK}.  There are many other results due to Carrell, Peterson,
and other authors related to (rational) smoothness of the Schubert
varieties.  For example, according to a result of D.~Peterson, 
see Carrell and Kuttler~\cite{CK}, smoothness of Schubert
varieties is equivalent to rational smoothness in the case of a
simply-laced root system.  Nevertheless none of these general criteria
give a simple efficient nonrecursive method (such as the
Lakshmibai-Sandhya criterion) for determining if a given Schubert
variety is smooth or not.  Recently, Billey~\cite{Bil} presented
analogues of Lakshmibai-Sandhya's theorem, for all classical types
$B_n$, $C_n$, and $D_n$.  However, these constructions, including the
definitions of patterns,  depend on a particular way to
represent elements in classical Weyl groups as signed permutations.

The main goals of this paper are to present a uniform approach to
pattern avoidance in general terms of root systems and to extend the
Lakshmibai-Sandhya criterion to the case of an arbitrary semisimple
Lie group $G$.  This approach using root subsystems will be
described in the next section.  Theorem~\ref{th:main} gives a
polynomial time algorithm for determining smoothness and rational
smoothness of Schubert varieties in $G/B$ in terms of root subsystems.
As a consequence of the main theorem, we get two additional criteria
for (rational) smoothness in terms of root systems embeddings and
double parabolic factorizations (see Theorems~\ref{th:main.embedding}
and~\ref{t:factored.form}).

Based on the ideas of root subsystems presented in this work, Braden
and the first author~\cite{BiBr} refined this notion and gave a lower
bound for the Kazhdan-Lusztig polynomials evaluated at $q=1$ in terms
of patterns.  They also introduce a geometrical construction which
identifies ``pattern Schubert varieties'' as torus fixed point
components inside larger Schubert varieties.  This can be used to give
another proof of one direction of our main theorem.  However, due to a
delay in publication, those results will appear first.

In Section~\ref{s:main}, we formulate our smoothness criterion
and describe the minimal lists of patterns needed to identify 
singular (rationally singular) Schubert varieties.
In Section~\ref{s:embeddings}, we present a computational
improvement using root system embeddings that reduces the minimal lists
to just 4 patterns (3 patterns) for (rational) smoothness test.  
The difference between smoothness
and rational smoothness is exhibited in the presence or absence of
rank 2 patterns. 
The connection to fully commutative elements is described in
Section~\ref{s:commutative}.  In Section~\ref{s:criteria}, we recall
several known characterizations of smoothness and rational smoothness
from the literature which we will use in the proof of the main
theorem.   In Section~\ref{s:lemmas}, we reformulate our main result
in terms of parabolic subgroups.  Then we prove two
statements on parabolic decomposition which will be used in the proof
of Theorem~\ref{th:main}, including Theorem~\ref{t:chains} which gives
a criterion for factoring Poincar\'e polynomials of Schubert varieties.
In Section~\ref{s:proof}, we give the details of the proof of the main
theorem.

\section{Root subsystems and the main results}\label{s:main}

As before, let $G$ be a semisimple simply-connected complex Lie group 
with a fixed Borel subgroup $B$.  Let $\mathfrak{h}$ be the Cartan subalgebra 
corresponding to a maximal torus contained in $B$.
Let $\Phi\in\mathfrak{h}^*$ be the corresponding root system,
and let $W=W_\Phi$ be its Weyl group.  The
choice of $B$ determines the subset $\Phi_+\subset \Phi$ of positive
roots.  The fact that a Schubert variety $X_w$, $w\in W$, in $G/B$ is
smooth (or rationally smooth) depends only on the pair $(\Phi_+,w)$.
We call such a pair {\it (rationally) smooth\/} whenever the
corresponding Schubert variety is (rationally) smooth.
The {\it inversion set\/} of an element $w$ in the Weyl
group $W_\Phi$ is defined by 
$$
I_\Phi(w) = \Phi_+\cap w (\Phi_-)\,,
$$
where $\Phi_-=\{-\alpha \mid \alpha\in \Phi_+\}$ is the set of negative roots.

The following properties of inversion sets are well-known,  
see~\cite[\S1, no~6]{Bou}.

\begin{lemma}
The inversion set $I_\Phi(w)$ uniquely determines the Weyl group
element $w\in W_\Phi$.  Furthermore, a subset $I\subseteq \Phi_+$ in
the set of positive roots is the inversion set $I_\Phi(w)$ for some
$w$ if and only if there exist a linear form $h$ on the vector space
$\mathfrak{h}^*$ such that $I=\{\alpha\in\Phi_+\mid h(\alpha)>0\}$.
\label{le:inversion-set}
\end{lemma}

A {\it root subsystem\/} of $\Phi$ is a subset of roots $\Delta\subset
\Phi$ which is equal to the intersection of $\Phi$ with a vector
subspace. Clearly, a root subsystem $\Delta$ is a root system 
itself in the subspace spanned by 
$\Delta$, see~\cite[\S1, no~1]{Bou}.
It comes with the natural choice of positive roots $\Delta_+=\Delta\cap
\Phi_+$.

By Lemma~\ref{le:inversion-set}, for any $w\in W_\Phi$ and any root 
subsystem $\Delta\subset\Phi$, the set of roots
$I_\Phi(w)\cap \Delta$ is the inversion set $I_\Delta(\sigma)$
for a unique element $\sigma\in W_\Delta$ in the Weyl group
of $\Delta$.
Let us define the {\it flattening map} $f_\Delta:W_\Phi\to W_\Delta$
by setting $f_\Delta(w)=\sigma$ where $\sigma$ is determined by its inversion
set $I_\Delta(\sigma) = I_\Phi(w)\cap \Delta$.

Recall that a graph is called a {\it star\/} if it is connected and it
contains a vertex incident with all edges.  Let us say that a root
system $\Delta$ is {\it stellar\/} if its Dynkin diagram is a star and
$\Delta$ is not of type $A_1$ or $A_2$.  For example, $B_{3}$ is
stellar but $F_{4}$ is not. Our first analogue of the
Lakshmibai-Sandhya criterion can be formulated as follows.  See also
Theorems~\ref{t:simply-laced}, \ref{t:21.patterns},
\ref{th:main.embedding} and~\ref{t:factored.form}.

\begin{theorem}\label{th:main}
Let $G$ be any semisimple simply-connected Lie group, 
$B$ be any Borel subgroup, with
corresponding root system $\roots$ and Weyl group $W=W_{\roots}$.  For
$w\in W$, the Schubert variety $X_w \subset G/B$ is smooth (rationally
smooth) if and only if, for every stellar root subsystem $\Delta$ in
$\Phi$, the pair $(\Delta_+,f_\Delta(w))$ is smooth (rationally
smooth).
\end{theorem}

The proof of Theorem~\ref{th:main} appears in Section~\ref{s:proof}.
If $\Delta$ is a root subsystem in $\Phi$ and $\sigma=f_\Delta(w)$,
then we say that the element $w$ in $W_\Phi$ {\it contains the
pattern\/} $(\Delta_+,\sigma)$.  
It follows from Theorem~\ref{th:main} that an element in $W_\Phi$ 
containing a non-smooth (non-rationally-smooth) pattern is also 
non-smooth (non-rationally-smooth).
Another explanation of this fact for rational smoothness
based on intersection homology can be found in the  work of
Billey and Braden~\cite{BiBr} mentioned above.

Let us say that an element $w$ {\it avoids the pattern} 
$(\Delta_+,\sigma)$ if $w$ does not contain a pattern isomorphic 
to $(\Delta_+,\sigma)$.
Clearly, Theorem~\ref{th:main} implies that the set of (rationally) 
smooth elements $w\in W_\Phi$ can be described as the set
of all elements $w$ that avoid patterns of several types.  Since there are
finitely many types of stellar root systems, the list of forbidden
patterns is also finite.

\begin{figure}[ht]
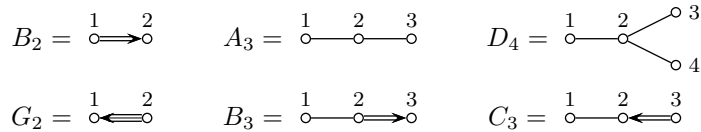

\pspicture(-100,-25)(150,22)
\rput(-100,10){$B_2=$}
\cnode(-80,10){2}{a1}
\rput(-80,17){{\fns 1}}
\cnode(-60,10){2}{a2}
\rput(-60,17){{\fns 2}}
\ncline[doubleline=true,doublesep=1pt]{->}{a1}{a2}

\rput(-20,10){$A_3=$}
\cnode(0,10){2}{v1}
\rput(0,17){{\fns 1}}
\cnode(20,10){2}{v2}
\rput(20,17){{\fns 2}}
\cnode(40,10){2}{v3}
\rput(40,17){{\fns 3}}
\ncline{-}{v1}{v2}
\ncline{-}{v2}{v3}

\rput(80,10){$D_4=$}
\cnode(100,10){2}{w1}
\rput(100,17){{\fns 1}}
\cnode(120,10){2}{w2}
\rput(120,17){{\fns 2}}
\cnode(140,0){2}{w3}
\rput(147,0){{\fns 4}}
\cnode(140,20){2}{w4}
\rput(147,20){{\fns 3}}
\ncline{-}{w1}{w2}
\ncline{-}{w3}{w2}
\ncline{-}{w4}{w2}

\rput(-100,-20){$G_2=$}
\cnode(-80,-20){2}{b1}
\rput(-80,-13){{\fns 1}}
\cnode(-60,-20){2}{b2}
\rput(-60,-13){{\fns 2}}
\ncline[doubleline=true,doublesep=2pt,linearc=.5]{<-}{b1}{b2}
\ncline[linewidth=.25pt]{-}{b1}{b2}

\rput(-20,-20){$B_3=$}
\cnode(0,-20){2}{B1}
\rput(0,-13){{\fns 1}}
\cnode(20,-20){2}{B2}
\rput(20,-13){{\fns 2}}
\cnode(40,-20){2}{B3}
\rput(40,-13){{\fns 3}}
\ncline{-}{B1}{B2}
\ncline[doubleline=true,doublesep=1pt,linearc=.5]{->}{B2}{B3}

\rput(80,-20){$C_3=$}
\cnode(100,-20){2}{C1}
\rput(100,-13){{\fns 1}}
\cnode(120,-20){2}{C2}
\rput(120,-13){{\fns 2}}
\cnode(140,-20){2}{C3}
\rput(140,-13){{\fns 3}}
\ncline{-}{C1}{C2}
\ncline[doubleline=true,doublesep=1pt,linearc=.5]{<-}{C2}{C3}

\endpspicture
\caption{Dynkin diagrams of stellar root systems}
\label{fig:stellar}
\end{figure}

Actually, the list of stellar root systems is relatively small:
$B_2$, $G_2$, $A_3$, $B_3$, $C_3$, and $D_4$. 
Figure~\ref{fig:stellar} shows their Dynkin diagrams labeled according
to standard conventions from~\cite{Bou}.
In order to use Theorem~\ref{th:main} as a (rational) smoothness test
we need to know all non-smooth and non-rationally-smooth 
elements in the Weyl groups with stellar root systems. The following table 
gives the numbers of such elements.

\bigskip
\begin{center}
\begin{tabular}{|lcccccc|}
\hline
stellar type:   &  $B_2$ & $G_2$ & $A_3$ & $B_3$ & $C_3$ & $D_4$ \\
\hline
non-smooth elements:    \mystrut
              &   1    & 5     & 2     & 20     &  20    & 49    \\ 
non-rationally-smooth elements: \mystrut
              &   0    & 0     & 2     & 14     &  14    & 49   \\[.05in]
\hline
\end{tabular}
\end{center}
\bigskip

There are several things to notice about the table.  In the simply-laced
cases $A_3$ and $D_4$ the numbers of non-smooth and non-rationally-smooth
elements coincide.
The rationally smooth elements in $B_n$ are exactly the same as the 
rationally smooth elements in $C_n$.  This explains why the number 14 
appears in both $B_3$ and $C_3$ cases.
Note that in general the number of non-smooth elements in $B_n$ is 
not equal to the number of non-smooth elements in $C_n$.  For examples,
we have 268 non-smooth elements in the  $B_4$ case and 270 non-smooth
elements in the $C_4$ case.

There are exactly two non-smooth elements of type $A_3$---they
correspond to the two forbidden patterns that appear 
Lakshmibai-Sandhya's criterion.  
Although there are 49 non-smooth elements of type $D_4$, 
only one (!)\ of these 49 elements contains no forbidden $A_3$ patterns. 
These three patterns (two of type $A_3$ and one of type $D_4$) 
are all patterns that are needed 
in the case of a simply-laced root system ($A$-$D$-$E$ case).

For all stellar types, let $s_1, s_2,\dots$ be the simple reflections
generating the corresponding Weyl groups labeled as shown on 
Figure~\ref{fig:stellar}.   Thus in both $A_3$ and $D_4$ cases
the reflection $s_2$ corresponds to the central node 
of the corresponding 
Dynkin diagram.  We will write elements of corresponding Weyl groups
as products of the generators~$s_i$.

\begin{theorem}\label{t:simply-laced}
Suppose that $\Phi$ is a simply-laced root system.  Then the Schubert variety 
$X_w$, $w\in W_\Phi$, is smooth if and only if $w$ avoids the following
three patterns:  two patterns of type $A_3$ given by the elements
$s_2 s_1 s_3 s_2$ and $ s_1 s_2 s_3 s_2 s_1$
and one pattern of type $D_4$ given by the element $s_2 s_1 s_3 s_4 s_2$.
\label{th:smooth-simply-laced}
\end{theorem}

Remark that, D.~Peterson has shown (unpublished, see Carrell 
and Kuttler~\cite{CK}) that in the simply-laced case a Schubert 
variety is smooth if and only if it is
rationally smooth.  Thus in the previous claim we can replace the word
``smooth'' by the phrase ``rationally smooth.''

For the case of arbitrary root systems (including non-simply-laced ones),
we need to list forbidden patterns of types $B_2$, $G_2$, $B_3$, and $C_3$.
The only non-smooth element of $B_2$ is $s_2 s_1 s_2$.  The non-smooth
elements of type $G_2$ are the 5 elements in the interval in the 
Bruhat order $[s_1s_2s_1,\wnot[$ ($\wnot$ is excluded), where $\wnot$
is the longest Weyl group element for type $G_2$.
There are also 6 non-smooth elements of type $B_3$ and 
6 non-smooth elements of type $C_3$ that contain no
forbidden $B_2$ patterns.  
The following theorem summarizes this data and gives the minimal 
list of patterns 
for the smoothness test.

We will write $[a,b,\dots,c]$ to denote the collection of words
$1,a,b,\dots,c$.  We concatenate this collection with another word as
follows: $[a,b,c]\,d$ is a shorthand for the four words $d, ad, bd,
cd$.

\begin{theorem} \label{t:21.patterns}
Let $\Phi$ be an arbitrary root system.  The Schubert variety 
$X_w$, $w\in W_\Phi$, is smooth if and only if $w$ avoids the  patterns
listed in the following table:
\smallskip
\begin{center}
\begin{tabular}{|c|c|c|}
\hline
stellar type   &  forbidden patterns & \# patterns \\
\hline
\hline
\mystrut
$B_2$  & $s_2s_1s_2$ & 1 \\[.05in]
\hline
\mystrut
$G_2$ & $[s_2]s_1s_2s_1[s_2]$,
      $s_1s_2s_1s_2s_1$ & 5
   \\[.05in]
\hline
\mystrut
$A_3$ & $s_2 s_1 s_3 s_2$, $s_1 s_2 s_3 s_2 s_1$ & 2 \\[.05in]
\hline
\mystrut
$B_3$ &
$s_2 s_1 s_3 s_2$, 
$s_1 s_2 s_3 s_2 s_1[s_3, s_3 s_2, s_2 s_3, s_2 s_3 s_2]$ & 6
\\[.05in]
\hline
\mystrut
$C_3$ & 
$[s_3]s_2 s_1 s_3 s_2[s_3]$,
$s_3 s_2 s_1 s_2 s_3$,
$s_1 s_2 s_3 s_2 s_1  s_3 s_2 s_3$
& 6 
\\[.05in]
\hline
\mystrut
$D_4$ &  $s_2 s_1 s_3 s_4 s_2$ & 1 \\[.05in]
\hline
\end{tabular}
\end{center}
\medskip
\label{th:smooth-general}
\end{theorem}

All Weyl group elements for the types $B_2$ and $G_2$ are rationally
smooth.  Thus we can ignore all root subsystems of these types in
rational smoothness test.  Rational smoothness can be defined in terms
of Kazhdan-Lusztig polynomials that depend only on the Weyl group.
The Weyl groups of types $B_3$ and $C_3$ are isomorphic.  Thus the
lists of non-rationally-smooth elements are identical in these two
cases.  The following theorem presents these lists.  

\begin{theorem} \label{t:min.patterns.rat.smooth}
Let $\Phi$ be an arbitrary root system.  The Schubert variety 
$X_w$, $w\in W_\Phi$, is rationally smooth if and only if 
$w$ avoids the patterns listed in the following table:
\smallskip
\begin{center}
\begin{tabular}{|c|c|c|}
\hline
stellar type   &  forbidden patterns & \# patterns\\
\hline
\hline
\mystrut
$A_3$ & $s_2 s_1 s_3 s_2$, $s_1 s_2 s_3 s_2 s_1$ & 2 \\[.05in]
\hline
\mystrut
$B_3=C_{3}$ & 
$[s_3]s_2 s_1 s_3 s_2 [s_3]$,
$[s_2] s_3 s_2 s_1 s_2 s_3 [s_2]$,
& 14
\\&
$s_1 s_2 s_3 s_2 s_1[s_3,s_2 s_3, s_3 s_2, s_2 s_3 s_2, s_3 s_2 s_3]$
&
\\[.05in]
\hline
\mystrut
$D_4$ & $s_2 s_1 s_3 s_4 s_2 $ & 1\\[.05in]
\hline
\end{tabular}
\end{center}
\medskip
\label{th:rational-smooth-general}
\end{theorem}

A quick glance on the tables in Theorems~\ref{th:smooth-general} 
and~\ref{th:rational-smooth-general} reveals that 
the lists of forbidden patterns for non-simply-laced cases
are longer than the list of three simply-laced forbidden
patterns.  In Section~\ref{s:embeddings} we show how to reduce the list 
patterns above to just the forbidden
patterns of types $B_2$, $A_3$, $D_{4}$ using embeddings of root systems.

\section{Root System Embeddings}
\label{s:embeddings}

In this section, we present an alternative notion of pattern avoidance
in terms of embedded root systems.  Again we can characterize
smoothness and rational smoothness of Schubert varieties.  The key
advantage of this approach is that we reduce the minimal number of
patterns to just 3 for rational smoothness and 4 for
smoothness.  While we believe this approach is useful for
computational purposes, we suspect root subsystems are better for
geometrical considerations.

Let $\Phi$ and $\Delta$ be two root systems in the vectors spaces $U$
and $V$, respectively.  An {\it embedding} of $\Delta$ into $\Phi$ is
a map $e:\Delta\to\Phi$ that extends to an injective linear map $U\to
V$.  
For example, any three positive roots $\alpha,\beta,$ $\gamma \in
\posroots$ define an $A_{3}$-embedding whenever $\alpha +\beta$,
$\beta+\gamma$ and $\alpha +\beta+\gamma$ are all in $\posroots$.  

Note that inner products are not necessarily preserved by embeddings as
they are with root subsystems.  Also note that every root subsystem
$\Delta$ in $\Phi$ gives an embedding of $\Delta$ into $\Phi$, but it
is not true that all embeddings come from root subsystems.  It is
possible that $\Delta$ embeds into $\Phi$ but the linear span of
$e(\Delta)$ in $\Phi$ contains some additional roots.  Nevertheless,
in the simply-laced case this can never happen.  For simply-laced root
systems, the notions of root subsystems and embeddings are essentially
equivalent.

We will say that a $k$-tuple of positive roots
$(\beta_1,\dots,\beta_k)$ in $\Phi$, gives a {\it $B_2$-embedding}, {\it
$A_3$-embedding}, or {\it $D_4$-embedding} if these vectors 
are the images of the simple roots in $\Delta$ 
for an embedding $\Delta\to\Phi$ with $\Delta$ 
of type $B_2$, $A_3$, or $D_4$, respectively.  For
example, $B_2$-embeddings are given by pairs of positive roots
$(\beta_1,\beta_2)$ such that both vectors $\beta_1+\beta_2$ and
$\beta_1+2\beta_2$ belong to $\Phi_+$.  Also {\it $A_3$-embeddings}
are given by triples of positive roots $(\beta_1,\beta_2,\beta_3)$
such that all vectors $\beta_1+\beta_2$, $\beta_2+\beta_3$, and
$\beta_1+\beta_2+\beta_3$ are roots in $\Phi_+$.

\begin{figure}[ht]
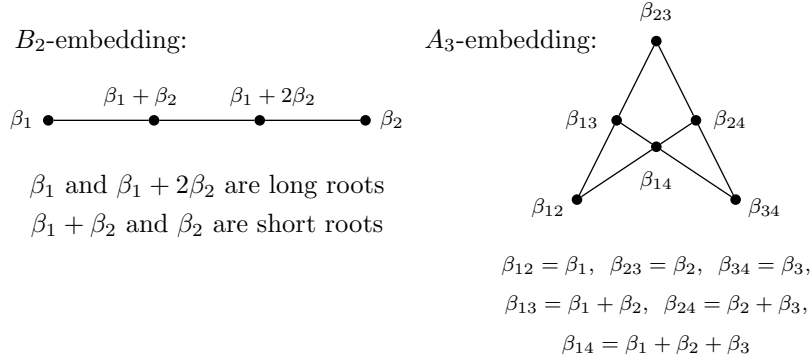

\pspicture(-200,-60)(60,75)
\rput(-180,60){$B_2$-embedding:}
\cnode*(-200,30){2}{a}
\rput(-210,30){\fns$\beta_1$}
\cnode*(-160,30){2}{b}
\rput(-165,39){\fns$\beta_1+\beta_2$}
\cnode*(-120,30){2}{c}
\rput(-115,39){\fns$\beta_1+2\beta_2$}
\cnode*(-80,30){2}{d}
\rput(-70,30){\fns$\beta_2$}
\ncline{-}{a}{b}
\ncline{-}{b}{c}
\ncline{-}{c}{d}
\rput(-140,5){$\beta_1$ and $\beta_1 + 2\beta_2$ are long roots}
\rput(-140,-10){$\beta_1+\beta_2$ and $\beta_2$ are short roots}

\rput(-25,60){$A_3$-embedding:}
\cnode*(0,0){2}{A}
\rput(-11,-3){\fns$\beta_{12}$}
\cnode*(30,60){2}{B}
\rput(30,71){\fns$\beta_{23}$}
\cnode*(60,0){2}{C}
\rput(71,-3){\fns$\beta_{34}$}
\cnode*(15,30){2}{AB}
\rput(02,30){\fns$\beta_{13}$}
\cnode*(45,30){2}{BC}
\rput(58,30){\fns$\beta_{24}$}
\cnode*(30,20){2}{ABC}
\rput(30,7){\fns$\beta_{14}$}
\ncline{-}{A}{AB}
\ncline{-}{AB}{B}
\ncline{-}{B}{BC}
\ncline{-}{BC}{C}
\ncline{-}{A}{ABC}
\ncline{-}{ABC}{BC}
\ncline{-}{C}{ABC}
\ncline{-}{ABC}{AB}
\rput(30,-25){\fns$\beta_{12}=\beta_1$, \ $\beta_{23}=\beta_2$, \ $\beta_{34}=\beta_3$,}
\rput(30,-40){\fns$\beta_{13}=\beta_1+\beta_2$, \ $\beta_{24}=\beta_2+\beta_3$,}
\rput(30,-55){\fns$\beta_{14}=\beta_1+\beta_2+\beta_3$}

\endpspicture
\caption{$B_2$- and $A_3$-embeddings}
\label{fig:A_3-embedding}
\end{figure}

Figure~\ref{fig:A_3-embedding} illustrates $B_2$- and
$A_3$-embeddings.  The vertices on the figure correspond to the
positive roots in the image of the embedding.  Here we used a
$(k-1)$-dimensional picture in order to represent collections of
$k$-dimensional vectors.  The vertices on the figure are the
intersections of the lines generated by the roots with a certain
affine hyperplane.  Therefore, inversion sets are determined by half
planes in these pictures.  A similar 3-dimensional figure can be
constructed for~$D_{4}$.  
Egon Schulte pointed out
\cite{schulte.personal} that the figure can be obtained by projecting
12 vertices of the regular 24-cell onto a tetrahedron spanned by 4
vertices of the 24-cell. In fact, it can be viewed as a model in
projective 3-space, and then it is actually related to the
half-24-cell.

The set of positive roots $\Phi_+$ in $\Phi$ and the embedding $e$
determines the set of positive roots $\Delta_+=e^{-1}(\Phi_+)$ in
$\Delta$.  We can extend the definition of the flattening map to
embeddings of root systems.  For an embedding $e:\Delta\to\Phi$, let
us define the {\it flattening map} $f_e:W_\Phi\to W_\Delta$ by setting
$f_e(w)=\sigma$, if the inversion set of $w$ pulls back to the
inversion set of $\sigma$ i.e.  $I_\Delta(\sigma) =
e^{-1}(I_\Phi(w))$.  According to Lemma~\ref{le:inversion-set}, the
element $\sigma$ is uniquely defined.

There are five $A_{3}$-embeddings into a root system $\Phi$ 
of type $B_{3}$ and seven into $\Phi$ of type $C_{3}$.
Among these twelve embeddings, five are necessary to classify
rationally singular Schubert varieties and three of these embeddings
lead to false positive classifications of rationally smooth elements
in $W_{C_{3}}$.  Therefore, we introduce the following definition
in order to eliminate the false conditions. 
For an embedding $e:\Delta\to\Phi$, let $\bar\Delta\subset\Phi$ be the
root subsystem in $\Phi$ spanned by the image $e(\Delta)$.
We say that an embedding $e:\Delta\to\Phi$ is {\it proper} if either
$\bar{\Delta}$ is not of type $B_{3},C_{3}$ or  
$\bar{\Delta}$ is of type $B_{3},C_{3}$  and 
there exists a
$B_{2}$-embedding $\epsilon :B_{2} \rightarrow \bar{\Delta}$ such that 
\begin{enumerate}
\item If $B_{2}$ has basis $\beta_{1},\beta_{2}$, then
$\epsilon(\beta_{1}+\beta_{2})=e(\alpha_{i})$ for some simple root $
\alpha_{i} \in \Delta$.  
\item  We have 
$\epsilon^{-1}(e(\Delta))=I_{B_{2}}(s_{2}s_{1}s_{2})=\{\beta_{1}+\beta_{2},\beta_{1}+2\beta_{2},
\beta_{2}\}$.   In words, the
image of the $B_{2}$-embedding intersects the image of the
$\Delta$-embedding in exactly three roots which correspond to the
inversion set of the unique singular Schubert variety
$X(s_{2}s_{1}s_{2})$ of type $B_{2}$.
\end{enumerate}
The root systems  $B_{3}$ and $C_{3}$ each have three 
$B_{2}$-embeddings and each of
these embeddings corresponds to exactly one proper $A_{3}$-embedding.

For an element $w$ in the Weyl group $W_\Phi$, we say that $w$
contains an {\it embedded pattern} of type $B_2$, $A_3$, or $D_4$ if
there is a proper embedding $e:\Delta\to\Phi$ such that

\begin{itemize}
\item
$\mathbf{B_2}$\,:  $\Delta$ is of type $B_{2}$ and 
$f_e(w) = s_2 s_1 s_2$;
\item$\mathbf{A_3}$\,: 
$\Delta$ is of type $A_{3}$ and 
$f_e(w) = s_1 s_2 s_3 s_2 s_1$ or 
$f_e(w) = s_2 s_1 s_3 s_2$;
\item$\mathbf{D_4}$\,:   
$\Delta$ is of type $D_{4}$ and 
and $f_e(w) = s_2 s_1 s_3 s_4 s_2$.
\end{itemize}
Recall here that the Coxeter generators $s_i$ of Weyl groups of types
$B_2$, $A_3$, and $D_4$ are labeled as shown on
Figure~\ref{fig:stellar}.  Note, that the reduced expressions above
are all of the form: central node conjugated by its neighbors or
neighbors conjugated by the central node.

Let $\Phi^\vee$ be the root system dual to $\Phi$.  
Its Weyl group $W_{\Phi^\vee}$ is naturally isomorphic to $W_\Phi$.
For an element $w\in W_\Phi$, we say that $w$ contains a 
\textit{dual embedded pattern} whenever the corresponding 
element in $W_{\Phi^\vee}\simeq W_\Phi$ contains an 
embedded pattern given by a proper embedding $e:\Delta\to\Phi^\vee$.

\begin{theorem}  \label{th:main.embedding}
Let $G$ be any semisimple simply-connected Lie group, $B$ be any Borel
subgroup, with corresponding root system $\roots$ and Weyl group
$W=W_{\roots}$. 
\begin{enumerate}
\item  For $w\in W$, the Schubert variety $X_w$ is
rationally smooth if and only if $w$ has no embedded patterns or dual
embedded patterns of 
types $A_3$ or $D_4$.  \label{th:embedded-patterns-rational}
\item For $w\in W_\Phi$, the Schubert variety $X_w$ is smooth
if and only if $w$ has neither embedded patterns 
of types $B_2$, $A_3$, or $D_4$, 
nor dual embedded patterns of types $A_3$ or $D_4$.
\label{th:embedded-patterns}
\end{enumerate}
\end{theorem}

Note that, the element $w$, corresponding to a smooth 
Schubert variety $X_w$, may contain 
dual embedded patterns of type $B_2$.
Thus smoothness of Schubert varieties, unlike rational smoothness, 
is not invariant with respect to duality of root system.

\begin{proof}
Any $B_{2}$, $A_{3}$, or $D_{4}$-embedding spans a root subsystem
whose rank must be at most 4.  Therefore, this theorem follows
directly from Theorem~\ref{th:main} by checking all root systems of
rank at most 4.
\end{proof}

\begin{figure}[ht]
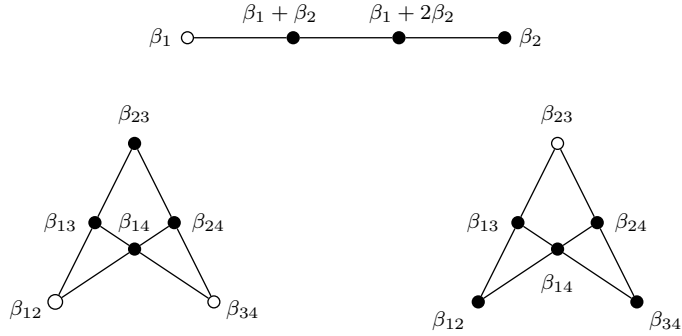

\pspicture(0,-35)(220,115)

\cnode(50,100){2.5}{aa}
\rput(40,100){\fns$\beta_1$}
\cnode*(90,100){2.5}{bb}
\rput(85,109){\fns$\beta_1+\beta_2$}
\cnode*(130,100){2.5}{cc}
\rput(135,109){\fns$\beta_1+2\beta_2$}
\cnode*(170,100){2.5}{dd}
\rput(180,100){\fns$\beta_2$}
\ncline{-}{aa}{bb}
\ncline{-}{bb}{cc}
\ncline{-}{cc}{dd}


\cnode(0,0){3}{A}
\rput(-11,-3){\fns$\beta_{12}$}
\cnode*(30,60){2.5}{B}
\rput(30,71){\fns$\beta_{23}$}
\cnode(60,0){2.5}{C}
\rput(71,-3){\fns$\beta_{34}$}
\cnode*(15,30){2.5}{AB}
\rput(02,30){\fns$\beta_{13}$}
\cnode*(45,30){2.5}{BC}
\rput(58,30){\fns$\beta_{24}$}
\cnode*(30,20){2.5}{ABC}
\rput(30,30){\fns$\beta_{14}$}
\ncline{-}{A}{AB}
\ncline{-}{AB}{B}
\ncline{-}{B}{BC}
\ncline{-}{BC}{C}
\ncline{-}{A}{ABC}
\ncline{-}{ABC}{BC}
\ncline{-}{C}{ABC}
\ncline{-}{ABC}{AB}

\cnode*(160,0){2.5}{a}
\rput(149,-7){\fns$\beta_{12}$}
\cnode(190,60){2.5}{b}
\rput(190,71){\fns$\beta_{23}$}
\cnode*(220,0){2.5}{c}
\rput(231,-7){\fns$\beta_{34}$}
\cnode*(175,30){2.5}{ab}
\rput(162,30){\fns$\beta_{13}$}
\cnode*(205,30){2.5}{bc}
\rput(218,30){\fns$\beta_{24}$}
\cnode*(190,20){2.5}{abc}
\rput(190,7){\fns$\beta_{14}$}
\ncline{-}{a}{ab}
\ncline{-}{ab}{b}
\ncline{-}{b}{bc}
\ncline{-}{bc}{c}
\ncline{-}{a}{abc}
\ncline{-}{abc}{bc}
\ncline{-}{c}{abc}
\ncline{-}{abc}{ab}
\endpspicture
\caption{Forbidden embedded patterns of types $B_2$ and $A_3$}
\label{fig:A_3-patterns}
\end{figure}

We mention one more computational simplification in applying
Theorem~\ref{th:main.embedding}.  For any $w\in W_\Phi$, there exists
a hyperplane that separates the sets of roots $I_\Phi(w)$ and
$\Phi_+\setminus I_\Phi(w)$. Figure~\ref{fig:A_3-patterns} illustrates
embedded patterns of types $B_2$ and  $A_3$. It is easy to see that each of these
inversion sets is determined by a half plane.  The black vertices
``$\,\pspicture(-2,-2)(2,2)\pscircle*(0,0){2.5}\endpspicture\,$''
correspond to the roots in the inversion set $I_\Phi(w)$ and the white
vertices
``$\,\pspicture(-2,-2)(2,2)\pscircle(0,0){2.5}\endpspicture\,$''
correspond to the roots outside the inversion set $I_\Phi(w)$.

Therefore, in order to search for embedded patterns for $w$ of types
$B_{2}$, $A_{3}$ and $D_{4}$ we only need to look for pairs, triple or
quadruples of the following forms:

\begin{enumerate}
\item[$\mathbf{B_2}$:] 
A pair of positive roots $(\beta_1,\beta_2)$ 
which forms the basis of a $B_2$-embedding
such that $\beta_1\not\in I_\Phi(w)$ and $\beta_1+\beta_2\in I_\Phi(w)$. 
\item[$\mathbf{A_3}$:] 
A triple of positive roots $(\beta_1,\beta_2,\beta_3)$ 
which forms the basis of a {\it proper\/} $A_3$-embedding 
such that (1) $\beta_{12},\, \beta_{34}\not\in I_\Phi(w)$ and $\beta_{14}\in I_\Phi(w)$;
or (2) $\beta_{23}\not\in I_\Phi(w)$ and $\beta_{13},\, \beta_{24}\in I_\Phi(w)$;
\item[$\mathbf{D_4}$:] 
A 4-tuple of positive roots $(\beta_1,\beta_2,\beta_3,\beta_4)$
which forms the basis of a  $D_4$-embedding such that 
$\beta_1+2\beta_2+\beta_3+\beta_4\in I_\Phi(w)$ and
$\beta_1+\beta_2+\beta_3,\ \beta_1+ \beta_2+\beta_4,\ 
\beta_2+\beta_3+\beta_4\not\in I_\Phi(w)$.
\end{enumerate}

\section{Other Elements Characterized by Pattern Avoidance}
\label{s:commutative}

In a series of papers (see~\cite{Fan}, \cite{Stemb1}, \cite{Stemb2} 
and reference wherein), Fan and Stembridge have developed a theory of
{\it fully commutative elements} in arbitrary Coxeter groups.  
By definition, an element in a Coxeter group is fully commutative
if all its reduced decompositions can be obtained from each other
by using only the Coxeter relations that involve commuting generators.

According to~\cite{BJS} the fully commutative elements in type $A$ 
are exactly the permutations avoiding the pattern $321$.  In
types $B$ and $D$, Stembridge has shown that the fully commutative
elements can again be characterized by pattern 
avoidance~\cite[Theorems 5.1 and 10.1]{Stemb1}

We note here that fully commutative elements are easily characterized
by root subsystems as well.  The following is an unpublished theorem
originally due to Stembridge~\cite{stem.personal}.  

\begin{proposition}\label{p:fully.commutative}
Let $W$ be any Weyl group with corresponding root system $\Phi$.  Then
$w \in W$ is fully commutative if and only if for every  root subsystem
$\Delta$ of type $A_{2}$, $B_{2}$, or $G_{2}$ we have 
$f_{\Delta}(w) \neq \wnot^{\Delta}$ where $\wnot^{\Delta}$ is the unique
longest element of $W_{\Delta}$.  In other words, $w$ is fully commutative
if and only if $w$ avoids the patterns given by the longest elements in 
rank $2$ irreducible root systems.
\end{proposition}

\begin{remark}
{\rm Fan, Stembridge and Kostant also investigated {\it abelian
elements} in Weyl groups.  An element $w\in W$ is abelian if its
inversion set $I(w)$ contain no three roots $\alpha,\beta$, and
$\alpha+\beta$.  Equivalently, $w\in W$ is abelian if the Lie algebra
$\mathfrak{b}\cap w(\mathfrak{b}_-)$ is abelian, where $\mathfrak{b}$
is Borel and $\mathfrak{b}_-$ is opposite Borel algebras.  For
simply-laced root systems, the set of abelian elements coincides with
the set of fully commutative elements.  The set of abelian elements
has a simple characterization in terms of embedded patterns.  Indeed,
by definition, $w\in W_\Phi$ is abelian if and only if there is no
$A_2$-embedding $e:\Delta\to\Phi$ such that the flattening $f_e(w)$ is
the longest element $\wnot^\Delta$ of~$W_\Delta$.  }
\end{remark}

\section{Criteria for Smoothness and Rational Smoothness}\label{s:criteria}

In this section, we summarize the three criteria for smoothness and
rational smoothness we rely on for the proof of Theorem~\ref{th:main}.

Let $\alpha_1,\dots,\alpha_n$ be the simple roots in $\Phi$ and 
let $\mathbb{Z}[\mathfrak{h}]$ denote the symmetric algebra generated
by $\alpha_1,\dots,\alpha_n$.  
For any $w,v\in W$ such that $w\leq v$, let us define 
$K_{w,v} \in \mathbb{Z}[\mathfrak{h}]$ by the recurrence
$$
\begin{array}{l}
\displaystyle
K_{w,w} = \prod_{\alpha \in \inv{w}} \alpha
\qquad\textrm{for $w=v$; }\\[.3in]
\displaystyle
K_{w,v}= 
\left\{
\begin{array}{cl}
K_{ws_{i},v} &	\textrm{ if }v< vs_{i}\\
K_{ws_{i},v} + (ws_{i}\alpha_{i})K_{ws_{i},vs_{i}}
&\textrm{ if } v > vs_{i}
\end{array}
\right.
\end{array}
$$
for $v\leq w$ and any simple reflection $s_{i}$ such that $ws_{i}<w$.
Then $K_{w,v}$ is a polynomial of degree $\ell(v)$ in the simple roots with non-negative
integer coefficients.  These polynomials first appeared in the work of
Kostant and Kumar~\cite{KK} on the nil Hecke ring, 
see~\cite{Bil2} for the recurrence.  

Kumar has given very general criteria for smoothness and rational
smoothness in terms of the nil-Hecke ring.  Through a series of
manipulations which were given in~\cite{BL}, one can obtain the
following statement from Kumar's theorem for finite Weyl groups.
Kumar's theorem in full generality applies to the Schubert varieties
for any Kac-Moody group.  However we would need to work with rational
functions of the roots.

\begin{theorem}\label{t:kumar}{\rm \cite{Kum,BL}}
Given any $v,w \in
W$ such that $v \leq w$, the Schubert variety $X_{v\wnot}$ is smooth
at $e_{w\wnot}$ if and only if
\begin{equation} \label{e:kumar.eq.1}
K_{w,v} = \prod_{\alpha
\in Z(w,v)} \alpha.  
\end{equation}   
where $Z(w,v) = \{\alpha \in \Phi_+ \, : \, v \not\leq
s_{\alpha} w \}$.
\end{theorem}

We can simplify the computations in Theorem~\ref{t:kumar} by
evaluating this identity at a well chosen point.  The modification
reduces the problem from checking a polynomial identity to checking
degrees plus a numerical identity.  Checking the degrees can be done
with a polynomial time algorithm since this only depends on the number
of positive roots.  

\begin{lemma}\label{l:kumar.reduction}
Let $r\in \mathfrak{h}$ be any regular dominant integral weight
i.e. $\alpha(r)\in \mathbb{N}_{+}$ for each $\alpha \in \Phi_{+}$.
Given any $v,w \in W$ such that $v \leq w$, the Schubert variety
$X_{v\wnot}$ is smooth at $e_{w\wnot}$ if and only if
\begin{equation} \label{e:billey-shimo} 
|Z(w,v)| = \ell(v) \text{\hspace{.2in} and \hspace{.2in} }
K_{w,v} (r) = \prod_{\alpha
\in Z(w,v)} \alpha (r).
\end{equation} \end{lemma}

\begin{proof}

We just need to prove the equivalence of \eqref{e:kumar.eq.1} and
\eqref{e:billey-shimo}.  Dyer~\cite{Dy} has shown that $K_{w,v}$ is
divisible by $\prod_{\alpha \in Z(w,v)} \alpha.$ Therefore,
$K_{w,v}=\prod_{\alpha \in Z(w,v)} \alpha$ if and only if their quotient
is 1.  We can check that the quotient is 1 by checking the degrees are
equal 
\[
\ell(v)=\mathrm{deg}(K_{w,v})
= \mathrm{deg}\left(\prod_{\alpha \in Z(w,v)} \alpha  \right)
=|Z(w,v)|, 
 \]
(in which case the quotient is a constant) and that $K_{w,v} (r)
= \prod_{\alpha \in Z(w,v)} \alpha (r)$.
\end{proof}

\begin{remark}
\label{r:roundoff.error}
Note, by choosing $r$ such that $\alpha(r)$ is always an integer, we do
not have to consider potential round off errors when checking equality.  
\end{remark}

\begin{remark}
\label{r:smooth.at.id} A Schubert variety $X_{w}$ is smooth at every
point if and only if it is smooth at $e_{\mathrm{id}}$.  Therefore, we
only need to check $\displaystyle K_{\wnot,w \wnot}(r)=\prod_{\alpha
\in Z(\wnot,w \wnot)} \alpha(r)$ when $|Z(\wnot,w\wnot)|=\ell(w\wnot)$ or
equivalently $|\{\alpha \in \posroots:\, s_{\alpha}\leq w \}| = \ell(w)$.
\end{remark}

The next criterion due to Carrell-Peterson is for rational smoothness.
The \textit{Bruhat graph $B(w)$ for $w \in W$} is the graph with vertices
$\{x \in W : x\leq w \}$ and edges between $x$ and $y$ if
$x=s_{\alpha}y$ for some $\alpha \in \posroots$ where $s_{\alpha}$ is
the reflection corresponding to $\alpha$,
\[
s_{\alpha} v = v - \frac{(v,\alpha)}{2(\alpha, \alpha)} \alpha.
\]
Note the Bruhat graph contains the Hasse diagram of the lower order ideal below $w$ in Bruhat order plus
some extra edges.  Let $$P_{w}(t) = \sum_{v \leq w} t^{\ell(v)},$$ then
$P_{w}(t^{2})$ is the \textit{Poincar\'e polynomial} for the
cohomology ring of the Schubert variety $X_w$.

\begin{theorem}\label{t:C-P}{\rm \cite{Car}}
The following are equivalent: 
\begin{enumerate}
\item $X_{w}$ is rationally smooth at every point.
\item The Poincar\'e polynomial $P_w(t)=\sum_{v\leq w}t^{\ell(v)}$ of
$X_{w}$ is symmetric (palindromic).
\item The Bruhat graph $B(w)$ is regular of degree $\ell(w)$, i.e., 
every vertex in $B(w)$ is incident to $\ell(w)$ edges.
\end{enumerate}
\end{theorem}

We can relate this theorem to the inversion sets $\inv{w}$ using the
following simple lemma, see~\cite{Bou}.

\begin{lemma}\label{l:inversion.sets}  
Fix a reduced expression $s_{a_{1}}s_{a_{2}}\dotsb s_{a_{p}}=w \in W$.
Let $\beta_{1},\dots, \beta_{n}$ be the simple roots in $\posroots$.
The following sets are all equal to the inversion set $\inv{w}$:
\begin{enumerate}
\item $\displaystyle \Phi_+\cap w (\Phi_-)$
\item $\displaystyle \{\alpha \in \Phi_+ : s_{\alpha} w <w \}.$
\item $\displaystyle \{s_{a_{1}}s_{a_{2}}\dotsb
        s_{a_{j-1}}\beta_{a_{j}}:\ 1\leq j\leq p \}$
\end{enumerate}
\end{lemma}


Let us label an edge $(x,s_\alpha x)$ in $B(w)$ by 
the root $\alpha\in\Phi_+$.
Then, by Lemma~\ref{l:inversion.sets}, the edges adjacent to $w$ in the
Bruhat graph $B(w)$ are labeled by the elements of $\inv{w}$ so
the degree $\mathrm{deg}(w)$ of the vertex $w$ is $\ell(w)$.
At any other vertex $x<w$ we know $\# \{\alpha \in \roots_{+}:
s_{\alpha }x < x \} = \ell(x)$ so $\mathrm{deg}(x)=\ell(x)+\# \{\alpha \in
\roots_{+}: x <s_{\alpha }x \leq w \}$.  Therefore, we have the
following lemma.

\begin{lemma}\label{l:bruhat.graph}
The Bruhat graph $B(w)$ is not regular 
if and only if there exists an $x<w$ such that
$$
\mathrm{deg}(x)>\mathrm{deg}(w) \iff \# \{\alpha \in \roots_{+} \mid 
x <s_{\alpha }x \leq w \} > \ell(w)-\ell(x).
$$
\end{lemma}

\section{Parabolic Decomposition}
\label{s:lemmas}

In the first lemma below, we give an alternative characterization of
pattern containment in terms of a parabolic factorization.  This leads
to an alternative characterization of smooth and rationally smooth
elements in the Weyl group.  We also give a method for factoring some
Poincar\'e polynomials of Schubert varieties.

Fix a subset $J$ of the simple roots.  Let $\Phi^{J}$ be the
root subsystem spanned by roots in $J$, and let
$\Phi^{J}_+=\Phi^J\cap \Phi_+$ be its set of positive roots.  
Let $W_{J}$ be the parabolic
subgroup generated by the simple reflections corresponding to $J$.
Let $W^{J}$ be the set of minimal length coset representatives for
$W_J\backslash W$ (moding out on the left). 
In other words,
\begin{equation}\label{e:min.reps}
W^J=\{v\in W\mid v^{-1}(\alpha)\in\Phi_+ \textrm{ for any } 
\alpha\in\Phi^J_+\}.
\end{equation}
Every $w \in W$ has a unique parabolic decomposition as the product
$uv=w$ where $u\in W_{J}$, $v\in W^{J}$ and $\ell(w)=\ell(u)+\ell(v)$, and
conversely, every product $u\in W_{J}$, $v\in W^{J}$ has
$\ell(uv)=\ell(u)+\ell(v)$ \cite[Prop.1.10]{Hum}.  Equivalently, if $w=uv$ is
the parabolic decomposition and $s_{a_{1}}s_{a_{2}}\dotsb s_{a_{p}}$,
$s_{b_{1}}s_{b_{2}}\dotsb s_{b_{q}}$ are reduced expressions for $u,v$
respectively then each $s_{a_{i}} \in W_{J}$ and
$s_{a_{1}}s_{a_{2}}\dotsb s_{a_{p}}s_{b_{1}}s_{b_{2}}\dotsb s_{b_{q}}$
is a reduced expression for $w$.

Let $\Delta \subset \roots$ be any root subsystem.  It was shown in
\cite{BiBr} that $\Delta$ is conjugate to $\Phi^{J}$ for some subset
$J$ of the simple roots, i.e. there exists a $v_1\in W^J$ such that
$v_1(\Delta)=\Phi^J$.  Clearly, $W_{\Delta}$ and $W_{J}$ are
isomorphic subgroups since the Dynkin diagrams for $\Delta$ and
$\Phi^{J}$ are isomorphic as graphs.  If there exist multiple
isomorphisms, any one will suffice.

\begin{lemma}\label{l:factored.form}
Let $\Delta \subset \roots$ be any root subsystem.  
Suppose that $v_1(\Delta)=\Phi^J$ for $v_1\in W^J$ as above.  
Let $w \in W$.  Let $u \in W_{J}$ and let $u'=v_1^{-1}u v_1$ 
be the corresponding
element in $W_\Delta$ under the natural isomorphism. 
Then $\fl_{\Delta} (w)=u'$ if
and only if there exits $v_{2} \in W^{J}$ such that $w=v_1^{-1}uv_2$.
\end{lemma}

\begin{proof}
The element $v_1 \in W^{J}$ gives a one-to-one correspondence 
$\alpha\mapsto
v_1(\alpha)$ between {\it positive} roots in $\Delta$ and {\it
positive} roots in $\Phi^J$.  Also for any $v_2\in W^J$, $v_2^{-1}$
maps positive roots of $\Phi^J$ to positive roots in $\Delta \subset
\Phi$ and negative roots of $\Phi^J$ to negative roots of $\Delta
\subset \Phi$.

Let $v_{2}=u^{-1}v_{1}w$.  We claim $v_{2} \in W^{J}$ since for any
$\alpha \in \Phi_{J}^{+}$, $v_{2}^{-1}(\alpha)=w^{-1}
v_{1}^{-1}u(\alpha) >0$.  Then we have a bijection from the inversions
of $u$ to the inversions of $w$ in $\Delta$:
\begin{align*}
\alpha\in\Delta_+ \cap \inv{w} 
& \iff w^{-1} (\alpha)   <0\\
& \iff v_{2}^{-1}u^{-1}v_{1} (\alpha)   <0\\
& \iff  u^{-1}v_{1}(\alpha ) \in \Phi_{-}^{J} \\
&\iff  v_{1}(\alpha) \in \inv{u}.
\end{align*}

\end{proof}

\begin{theorem}\label{t:factored.form}
Let $\Phi$ be any root system and let $J_1,\dots,J_s$ be a collection 
of subsets of simple roots such that all parabolic subsystems 
$\Phi^{J_1},\dots,\Phi^{J_s}$ 
are stellar and they include all possible stellar types present in the 
Dynkin diagram of $\Phi$. 
Then $w\in W$ is (rationally) smooth if and only if it cannot be
presented in the form $w=v_1^{-1}uv_2$, where  
$v_1,v_2\in W^{J_i}$, $u\in W_{J_i}$, for $i\in\{1,\dots,s\}$, 
and $u$ is (rationally) singular element in $W_{J_i}$. 
\label{thm:W-W_J}
\end{theorem}

\begin{proof}
Suppose that $\Delta$ is any stellar root subsystem in $\Phi$.  Then
$\Delta$ is conjugate to $\Phi^{J}$, where $J=J_i$ for some $i$.  Now
Lemma~\ref{l:factored.form} shows that Theorem~\ref{thm:W-W_J} is
equivalent to Theorem~\ref{th:main}.
\end{proof}

Each stellar parabolic subset $J_i$ in Theorem~\ref{t:factored.form}
consists of a node in the Dynkin diagram together with its neighbors.
We need to pick all nonisomorphic such subsets.
For example, $s=1$, for $\Phi$ of type $A_n$;  $s=2$, for any other
simply laced type; and $s=3$ for $\Phi$ of type $B_n$ or $C_n$ with $n\geq 4$.

Theorem~\ref{thm:W-W_J} implies the following statement.

\begin{corollary}
Let us fix any subset of simple roots $J$.  
Suppose that $u\in W_J$ and $v_1,v_2\in W^J$ are
such that $v_1^{-1}uv_2$ is a (rationally) smooth element in $W$.
Then $u$ is a (rationally) smooth in element in $W_J$.
\label{cor:W_J}
\end{corollary}

It was shown in~\cite{BFL} that,
for any $w\in W$ and a subset $J$ of the simple roots,
the parabolic subgroup $W_J$ has a unique maximal element 
$m(w,J)\in W_J$ below $w$ in the Bruhat order. 
The following theorem generalizes the factoring formulas for
Poincar\'e polynomials found in~\cite{Gash,Bil} and~\cite[Thm.11.23]{BL}.  
Using this theorem one can simplify the search
for the palindromic Poincar\'e polynomials which appear in
Theorem~\ref{t:C-P}.

\begin{theorem}\label{t:chains}
Let $J$ be any subset of the simple roots. Assume $w\in W $ has the
parabolic decomposition $w=u\cdot v$ with $u\in W_{J}$ and $v \in
W^{J}$ and furthermore, $u=m(w,J)$.  Then
$$
P_{w}(t)=P_{u}(t)\,P^{W^{J}}_{v}(t) 
$$
where $\displaystyle P^{W^{J}}_{v}(t) =
\sum_{z \in W^{J},\, z\leq v } t^{\ell(z)}$ 
is the Poincar\'e polynomial for $v$ in the quotient.  
\end{theorem}

\begin{proof}
Let $B(w)=\{x\in W \mid x \leq w \}$ and 
$B_{W^{J}}(v)=\{z\in W^{J} \mid z \leq v \}$.  
We will show there exists a rank
preserving bijection $f:B(w)\rightarrow B(u) \times B_{W^{J}}(v)$.  

 Given any $x \leq w$, say $x$ has
parabolic decomposition $x=yz$ with respect to $J$, then $y \leq
u=m(w,J)$ since $m(w,J)$ is the unique maximal element below $w$ and
in $W_{J}$. Furthermore, Proctor~\cite[Lemma 3.2]{pr} has shown that if
$z \in W$ is also a minimal length element in the coset $W_{J}z$, then
$z \leq w$ if and only if $z \leq v$.  Therefore, we can define a map
\begin{equation}\label{e:chains}
f : B(w) \rightarrow B(u) \times B_{W^{J}}(v)
\end{equation}
by mapping $x$ to $(y,z)$.  Note that this map is injective and rank
preserving since $\ell(x)=\ell(y)+\ell(z)$ by the properties of the parabolic
decomposition.

Conversely, given any $y \in W$ such that $y \leq u$ and given any $z
\in B_{W^{J}}(v)$, then actually $y \in W_{J} $ and we have $ yz \leq
w$ with $\ell(yz)=\ell(y)+\ell(z)$.  Therefore, $yz$ can be written as a
subexpression of the reduced expression for $w$ which is the
concatenation of reduced expressions for $u$ and $v$.  Furthermore,
$f(yz)=(y,z)$ since $z$ is the unique minimal length coset
representative in the coset containing $yz$.  Hence $f$ is surjective.
\end{proof}

Gasharov~\cite{Gash} and 
Billey~\cite{Bil} have shown that, for the classical types, 
the Poincar\'e polynomials of rationally smooth Schubert varieties 
have very nice factorizations.  Gasharov and Reiner~\cite{GR2000}, 
Ryan~\cite{ryan}, and Wolper~\cite{Wolper}
have in fact shown that smooth Schubert varieties can be described as
iterated fiber bundles over Grassmannians.  We see a similar phenomena
for the exceptional types.  


\begin{corollary}
Every Poincar\'e polynomial of a rationally smooth Schubert variety in
types $A_{n}, B_{n}, C_{n}, D_{n}, E_{6}, E_{7},
E_{8}, G_{2},F_{4}$ factors
into a product of symmetric factors each of which are Poincar\'e
polynomials indexed by elements in a maximal parabolic quotient
$W/W_{J}$.
\end{corollary}

We only use the following lemma in the proof of Theorem~\ref{th:main}
while proving our criteria for rational smoothness.

\begin{lemma}\label{l:reduction.2}
Let $w=uv$ be the parabolic decomposition of $w \in W$ with $u\in
W_{J}$ and $v\in W^{J}$.
\begin{enumerate}
\item If $X_u$ is not rationally smooth then $X_w$ is not
rationally smooth.
\item For any root subsystem $\Delta \subset \roots^{J}$, we have
$\inv{u}\cap \Delta = \inv{w}\cap \Delta$.   Therefore, if $u$ has a
pattern $(\Delta_+,\sigma)$ then so does $w$.
\end{enumerate}
\end{lemma}

Note, it is not true that if $X_v$ is singular in the quotient
$G/P_{J}$ then $X_w$ is necessarily singular in $G/B$.  Also, it is
possible for both $X_u$ to be smooth in $G/B$, $X_v$ smooth in the
quotient and yet $X_w$ to be singular in $G/B$ if $u \neq m(w,J)$.

\begin{proof}
To prove the first statement, assume $X_u$ is not rationally
smooth. Then there exits a vertex $x<u$ in the Bruhat graph for $u$
where $\mathrm{deg}(x)$ is too large by Theorem~\ref{t:C-P} and
Lemma~\ref{l:bruhat.graph}, namely
$$
\# \{\alpha \in \roots_{+}: x < s_{\alpha }x \leq u \} >
\ell(u)-\ell(x).
$$

We claim $\mathrm{deg}(xv)$ is too large for $X_w$.  Note, $x \leq
s_{\alpha}x \leq u$ implies $s_{\alpha} \in W_{J}$, and therefore $xv \leq
s_{\alpha }xv \leq w$ by the properties of the parabolic
decomposition.  Therefore,
$$
\begin{array}{l}
\# \{\alpha  \in \posroots: xv < s_{\alpha }xv \leq w \} 
\ \geq \  \# \{\alpha  \in \roots_{+}: x < s_{\alpha }x \leq u \}\\
>\ \ell(u)-\ell(x)
=\  \ell(w)-\ell(v)-\ell(x)
=\  \ell(w)-\ell(xv).
\end{array}
$$
Hence, $X_w$ is not rationally smooth, proving the first statement.

The second statement follows directly from
Lemma~\ref{l:factored.form}.  
\end{proof}

\section{Proof of Theorem~\ref{th:main}}
\label{s:proof}

The semisimple Lie groups come in 4 series: $A_{n}, B_{n}, C_{n},
D_{n}$ and 5 exceptional types $E_{6}, E_{7}, E_{8}, F_{4}, G_{2}$.
The proof of Theorem~\ref{th:main} for the four infinite series
follows easily on the known characterization for smooth Schubert
varieties in terms of pattern avoidance.  The proof for the
exceptional types was much more difficult.  One might imagine that a
routine verification would suffice for these finite Weyl groups.
However, both verifying smoothness in $F_{4}$ and rational smoothness
in $E_{8}$ directly would be impossible in our life time using
previously known techniques.  The exceptional types were proved with
the aide of a large parallel computer after making several reductions
in complexity.  These reductions in complexity also give insight into
the intricate geometry of the exceptional types.  

Recall, classical pattern avoidance is defined in terms of the
following function which {\it flattens} any subsequence into a signed
permutation.  Let $\mathcal{B}_{n}$ denote the signed permutation
group.  Elements in $\mathcal{B}_{n}$ can be written in
\textit{one-line notation} as an ordered list of the numbers
$1,\dots, n$ with a choice of sign for each entry.  For example, $3
\bar{2}1 = s_{2}s_{3}s_{2}s_{1}s_{2} \in \mathcal{B}_{n}$ (barred
numbers are negative).  The group $\mathcal{B}_{n}$ is 
isomorphic to the Weyl group of type $B_{n}$ and $C_{n}$.  
The Weyl group of type $D_{n}$ is the subgroup of $\mathcal{B}_{n}$
in which all elements have an even number of negative entries; and
the Weyl group of type $A_{n-1}$ (the symmetric group $S_{n}$) 
is the subgroup in which all elements have no negative entries.

\begin{definition}
Given any sequence $a_1 a_2 \ldots a_k$ of distinct non-zero real
numbers, define $\fl(a_1 a_{2} \ldots a_{k})$ to be the unique element
$b=b_1 \ldots b_k$ in $\mathcal{B}_k$ such that
\begin{itemize}
\item 
both $a_j$ and $b_j$ have the same sign.
\item for all $i,j$, we have $|b_i|<|b_j|$ if and only if $|a_i|<|a_j|$.
\end{itemize}
\end{definition}
For example, $\fl(\bar{6}, 3, \bar{7}, 1)=\bar{3}2 \bar{4} 1$.  Any
sequence containing the subsequence $\bar{6}, 3, \bar{7}, 1$ does not
avoid the pattern $\bar{3} 2 \bar{4} 1$.

\begin{theorem}\label{t:classical.patterns}{\rm \cite{Bil,LSa}}
Let $W$ be one of the groups $W_{A_{n-1}}, W_{B_{n}}, W_{C_{n}}$ or
$W_{D_{n}}$ and let $w \in W$.  
Then $X_{w}$ is (rationally) smooth if and
only if for each subsequence $1\leq i_1<i_2<i_3 <i_4 \leq n$,
$\fl(w_{i_1}w_{i_2}w_{i_3}w_{i_4})$ corresponds to a (rationally)
smooth Schubert variety.
\end{theorem}

In order to prove Theorem~\ref{th:main}, we claim
$\fl(w_{i_1}w_{i_2}w_{i_3}w_{i_4})=v$ if and only if
$I_{\subroots}(v)= I_{\roots}(w)\cap \subroots$ where $\subroots$ is
the root subsystem of type $B_{4}$ in the span of
$e_{|w_{i_1}|},e_{|w_{i_2}|},e_{|w_{i_3}|},e_{|w_{i_4}|}$.  We will
prove this claim in type $B$, the remaining cases are similar.  Then
verification of types $A_{3}$, $B_{4}$, $C_{4}$ and $D_{4}$
suffices to check the theorem in the classical case.

For type $B_n$, let us pick the linear basis 
$e_1,\dots,e_n$ in $\mathfrak{h}^*$ 
such that the simple roots are given by 
$e_{1}, e_{2}-e_{1},\dots, e_{n}-e_{n-1}$.  Then 
$\Phi_+ = \{e_{k}\pm e_{j}\ \mid 1 \leq j < k <
n \} \cup \{e_{j}\ :\ 1\leq j\leq n \}$.  A signed permutation $w$
acts on $\mathbb{R}^{n}$ by 
$$
w(e_{j})=\left\{
\begin{array}{cl}
e_{w_{j}} &\textrm{ if } w_{j}>0\\
-e_{|w_{j}|} &	\textrm{ if } w_{j}<0.  
\end{array}
\right.
$$

Explicitly, $\inv{w}=\posroots \cap w \negroots$ is the union of the
following three sets 
$$
\begin{array}{l}
  \{w(-e_{j}):\ w_{j}<0 \} \\
  \{w(e_{j}-e_{k}):\ j<k,\ w_{j}>|w_{k}| \} \\
 \{w(\pm e_{j}-e_{k}):\ j<k,\ w_{k}<0 \textrm{ and } |w_{j}|<|w_{k}| \}.
\end{array}
$$
Therefore, deciding if $w(-e_{j})$ or $w(\pm e_{j}- e_{k}) \in
\inv{w}$ depends only on the relative order and sign patterns on
$w_{j}$ and $w_{k}$.  By definition of the classical flattening
function $\fl(w_{i}w_{i_2}w_{i_3}w_{i_4})=v \in W_{B_{4}}$ if
$w_{i}w_{i_2}w_{i_3}w_{i_4}$ and $v_{1}v_{2}v_{3}v_{4}$ have the same
relative order and sign pattern.  Hence, when $\Delta$ is the root
subsystem of type $B_{4}$ determined by
$e_{|w_{i_1}|},e_{|w_{i_2}|},e_{|w_{i_3}|},e_{|w_{i_4}|}$, we have
$I_{\roots}(w)\cap \subroots = I_{\subroots}(v)$ if and only if
$\fl(w_{i_1}w_{i_2}w_{i_3}w_{i_4})=v$.  This proves the claim and 
finishes the proof of Theorem~\ref{th:main} for the classical types.

\medskip

Next consider the root systems of types $G_{2}$ and $F_{4}$.  We can
simply check Theorem~\ref{th:main} by computer using the modified
version of Kumar's criterion for determining smoothness and the
Carrell-Peterson criteria discussed in Section~\ref{s:criteria}.  In
particular, for $G_{2}$, we use Remark~\ref{r:smooth.at.id} to find
all singular elements Schubert varieties.  They are $X_{s_1 s_2s_1}$,
$X_{s_1 s_2s_1 s_2}$, $X_{s_2 s_1s_2 s_1}$, $X_{s_1 s_2s_1 s_2s_1}$,
$X_{s_2 s_1s_2 s_1 s_2}$ (assuming $\alpha_{1}$ is the short simple
root).  Pattern avoidance using root subsystems does not offer any
simplification of this list.  However, using root systems embeddings
in Section~\ref{s:embeddings}, these singular Schubert varieties all 
follow from one $B_{2}$ pattern.  All Schubert varieties of type
$G_{2}$ are rationally smooth.

For $F_{4}$ we use the following algorithm to verify
Theorem~\ref{th:main}:

\begin{enumerate}
\item Make a matrix of size $|W|^{2}$ containing the values
$K_{w,v}(r)$ computed recursively using the formula
$$
K_{w,v}(r)= \left\{
\begin{array}{cl}
K_{ws_{i},v}(r) &\textrm{ if }v< vs_{i}\\
K_{ws_{i},v}(r) + \left(ws_{i}\alpha_{i}(r)\right)K_{ws_{i},vs_{i}}(r)&
\textrm{ if }v> vs_{i}
\end{array}
\right.
$$
where $s_{i}$ is any simple reflections such that $ws_{i}<w$.
\item 
Identify all subsets $\{\gamma_{1},\dots , \gamma_{p} \} \subset
\posroots$ for $p=2,3,4$ such that $\{\gamma_{1},\dots , \gamma_{p}
\}$ forms a basis for a root subsystem $\Delta$ of type $B_{2}, A_{3},B_{3}$
or $C_{3}$ (no root subsystems of types $G_{2}$ or $D_{4}$ appear in
$F_{4}$).  Let $\bases$ be the list of all such root subsystem bases.  

\item For each such root subsystem $\Delta$ with basis $B \in \bases$, find all
singular  Schubert varieties $X_v$ using
Remark~\ref{r:smooth.at.id} and add $I_{\Delta}(v)$ to a list
called $\bad(B)$.  Note this list can be significantly simplified by
removing all $B_{3}, C_{3}$ Schubert varieties which are
classified as singular using a $B_{2}$ root subsystem.

\item\label{i:4} For each $w \in W_{F_{4}}$, check $\ell(w) = |\{\gamma: \, s_{\gamma }\leq w \}|$ and  $K_{\wnot,w
\wnot}(r)=\prod_{\gamma \in Z(\wnot,w \wnot)} \gamma(r)$ if and only
if no $B \in \bases$ exists such that $\inv{w}\cap \mathrm{span}(B)$
is a member of $\bad(B)$.
\end{enumerate}

To verify the theorem for rational smoothness in $F_{4}$ we use the
same algorithm except $\bad(B)$ should contain all of the inversion
sets for rationally singular Schubert varieties of types $A_{3},
B_{3}, C_{3}, D_{4}$ and in Step~\ref{i:4} use the palindromic
Poincar\'e polynomial criterion for rational smoothness from
Theorem~\ref{t:C-P}.

\medskip

Finally, for $E_{6}, E_{7}$ and $E_{8}$ it suffices to check the
theorem on $E_{8}$ since the corresponding root systems and Weyl
groups are ordered by containment.  Note, the Weyl group of type
$E_{8}$ has 696,729,600 elements so creating the matrix as in Step~1
above is not possible with the current technology.  Therefore, a
different method for verifying the main theorem was necessary.
Recall, D.~Peterson has shown, see~\cite{CK}, that smoothness and rational
smoothness are equivalent for simply laced Lie groups, i.e. $A_{n},
D_{n}$ and $E_{6},E_{7},E_{8}$.  Unfortunately, computing $P_{w}(t)$
for all $w \in E_{8}$ and applying either of the Carrell-Peterson
criteria in Theorem~\ref{t:C-P} is out of the question, however we
made the following observations:

\begin{enumerate}
\item In $E_{8}$, if $P_{w}(t)$ is not symmetric then approximately
99.989\% of the time the only coefficients we need to check are of
$t^{1}$ and $t^{\ell(w)-1}$.  In fact, in all of $E_{8}$, all one ever
needs to check is the first 6 coefficients (starting at $t^{1}$)
equals the last 6 coefficients.  Note, the coefficient of $t^{1}$ is
just the number of distinct generators in any reduced expression for
$w$ and the coefficient of $t^{\ell(w)-1}$ is the number of 
$v=ws_{\alpha}$ for $\alpha \in \Phi_+$ 
such that $\ell(v)=\ell(w)-1$ which can be efficiently computed.

\item In $E_{8}$, if $P_{w}(t)$ is symmetric then there always exists
a factorization according to Theorem~\ref{t:chains} where $J$ is a
subset of all simple roots except one which corresponds to a leaf of
the Dynkin diagram.  This factored formula makes it easy to check the
palindromic property recursively.

\item Let $J$ be the set of all simple roots in $E_{8}$ except
$\alpha_{1}$.   Here we are labeling the simple roots according to the
following Dynkin diagram:
\begin{center}
\pspicture(-20,-15)(120,27)
\rput(-20,10){$E_8=$}
\cnode(0,10){2}{v1}
\rput(0,17){{\fns 1}}
\cnode(20,10){2}{v3}
\rput(20,17){{\fns 3}}
\cnode(40,10){2}{v4}
\rput(40,17){{\fns 4}}
\cnode(40,-10){2}{v2}
\rput(33,-10){{\fns 2}}
\cnode(60,10){2}{v5}
\rput(60,17){{\fns 5}}
\cnode(80,10){2}{v6}
\rput(80,17){{\fns 6}}
\cnode(100,10){2}{v7}
\rput(100,17){{\fns 7}}
\cnode(120,10){2}{v8}
\rput(120,17){{\fns 8}}
\ncline{-}{v1}{v3}
\ncline{-}{v3}{v4}
\ncline{-}{v2}{v4}
\ncline{-}{v4}{v5}
\ncline{-}{v5}{v6}
\ncline{-}{v5}{v6}
\ncline{-}{v6}{v7}
\ncline{-}{v7}{v8}
\endpspicture
\end{center}
By Lemma~\ref{l:reduction.2}, we only need to test
$w=uv$ where $v \in W/W_{J}$, $u \in W_{J} \approx W_{D_{7}}$ and
$X_u$ smooth.  There are 9479 elements of $W_{D_{7}}$ which
correspond to smooth Schubert varieties and 2160 elements in
$W/W_{J}$. 
\end{enumerate}

\medskip

With these three observations, we can complete the verification of
$E_{8}$ using the following algorithm:

\begin{enumerate}
\item Identify all root subsystems of types $A_{3},D_{4}$ and their
bases as above (since $E_{8}$ is simply laced it has no root
subsystems of the other types).  
Call the list of bases $\bases$.

\item For each such root subsystem $\subroots$ with basis $B \in \bases$, find all
singular (or equivalently rationally singular) Schubert varieties
$X_v$ using Remark~\ref{r:smooth.at.id}  and add $I_{\subroots}(v)$
to a list called $\bad(B)$.  Note this list can be significantly
shortened by removing all $D_{4}$ Schubert varieties which are
classified as singular using an $A_{3}$ root subsystem. 

\item Identify all 2160 minimal length coset representatives in the
quotient $W_{J} \backslash W$ (moding out on the left).  Call this set
$Q$.  

\item Identify all 9479 elements in $\textrm{Smooth-}W_{J} \approx
\textrm{Smooth-}D_{7} := \{u\in W_{D_{7}}: X_u \textrm{ is smooth}\} $ 
using classical pattern avoidance or root subsystems.

\item For each $u \in \textrm{Smooth-}W_{J}$ and $v \in Q$, let $w=uv$
and check if  $B \in \bases$ exists such that 
$\inv{w}\cap \mathrm{span}(B) \in \bad(B)$. 
\begin{enumerate}
\item If yes, check as many coefficients as necessary to show
$P_{w}(t)$ is not symmetric.  Here 5 coefficients sufficed for all
$u,v$.
\item If no, attempt to factor $P_{w}(t)$ by taking $J'$ to be all
simple roots except one of the leaf nodes of the Dynkin diagram and if
$w=m(w,J')\,v'$ or $w^{-1}=m(w^{-1},J')\,v'$ in the corresponding
parabolic decomposition then apply Theorem~\ref{t:chains}.  For every
$w$ in this case, there exists some such $J'$ so that
$P_{w}(t)=P_{u'}(t)\,P_{v'}^{W^{J'}}(t)$ where $P^{W^{J'}}_{v'}(t)$ is
symmetric and $P_{u'}(t)$ factored into symmetric factors recursively
by peeling off one leaf node of the Dynkin diagram at a time.
\end{enumerate}

\end{enumerate}

This completes the proof of Theorem~\ref{th:main}.
Theorems~\ref{t:simply-laced}, \ref{t:21.patterns}
and~\ref{t:min.patterns.rat.smooth} follow directly.

\begin{conjecture}\label{c:checks} Let $\roots$ be a simply laced root
system of rank $n$.  Say $(\posroots,w)$ is not a rationally smooth
pair.  We conjecture that one only need to compare the first $n$
coefficients and the last $n$ coefficients of $P_{w}(t)$ in order to
find an asymmetry.  Equivalently, the Kazhdan-Lusztig polynomial
$P_{\mathrm{id},w}(q)$ has a non-zero coefficient among the terms
$q^{1}, q^{2},\dots, q^{n}$.
\end{conjecture}

We can show that for $A_{n}$ one only needs to check $n-2$
coefficients, for $D_{5}$ and $E_{6}$ one needs to check 3
coefficients, and for $E_{8}$ one needs to check 5 coefficients.  
For $F_4$ (which is not simply laced) one needs to check 3 coefficients,
and for $B_{5}$ one needs to check 6 coefficients.

\section{Acknowledgments} The authors would like to thank Bert
Kostant, V.~Lakshmibai, Egon Schulte and John Stembridge for their
insightful comments.  Also, we are very grateful to Paul Viola and
Mitsubishi Electric Research Labs for the use of their parallel
computer which ran our proof for $E_{8}$ and to Shrawan Kumar for
pointing out an false statement in an earlier draft.

As with Appel and Haken's 1976 proof of the Four Color Theorem, our
proof has met some criticism as to the merit of using a computer
verification.  In fact, it has been rejected by two journals based not
on the importance or originality of the results, but on the method of
proof.  We believe quite to the contrary that every significant
computer aided proof is a major accomplishment in expanding the role of
computers in mathematics.  It is like practicing to use induction
2000 years ago; it was a highly creative and influential achievement.
With years of practice, we have become quite proficient with the
induction technique.  However, computer aided proof is a fledgling
technique that certainly will have a major impact on the future of
mathematics.  Therefore, we hope our method of proof will actually
make a much broader impact on the future of mathematics than our
main theorem.

Perhaps the only thing a complete human proof could add is an
intuitive explanation for why stellar root subsystems contain all the
bad patterns.  This remains an open problem.

\end{document}